\documentclass[10pt]{article}

\usepackage{amssymb}

\usepackage{latexsym}

\usepackage{amsbsy}

\usepackage[all]{xypic}

\usepackage{amscd}

\title{Applications of BGP-reflection functors: isomorphisms of cluster algebras
\thanks{Supported by the NSF of China (Grants 10471071) and  by
the Cultivation Fund of the Key Scientific and Technical Innovation
Project, Ministry of Education of China.}}

\author{Bin Zhu\thanks{E-mail: bzhu@math.tsinghua.edu.cn}}

\date{ \small {Department of Mathematical Sciences, \\  Tsinghua University,
100084 Beijing, P. R. China}}

\begin{document}

\maketitle

\def\s{\stackrel}
\def\gama{\gamma}
\def\Longrightarrow{{\longrightarrow}}
\def\P{{\cal P}}
\def\A{{\cal A}}
\def\F{\mathcal{F}}
\def\X{\mathcal{X}}
\def\T{\mathcal{T}}
\def\m{\textbf{ M}}
\def\t{{\tau }}
\def\b{\textbf{d}}
\def\K{{\cal K}}

\def\G{{\Gamma}}
\def\e{\mbox{exp}}

\def\righta{\rightarrow}

\def\s{\stackrel}

\def\ncong{\not\cong}

\def\mathbb{\NN}

\def\Hom{\mbox{Hom}}
\def\Ext{\mbox{Ext}}
\def\ind{\mbox{ind}}
\def\coprod{\amalg }
\def\L{\Lambda}
\def\c{\circ}
\def\mu{\multiput}

\newcommand{\uHom}{\operatorname{\underline{Hom}}\nolimits}
\newcommand{\End}{\operatorname{End}\nolimits}
\renewcommand{\r}{\operatorname{\underline{r}}\nolimits}
\def \text{\mbox}

%Egne definisjoner%

\begin{center}

\begin{minipage}{12cm}{\footnotesize\textbf{Abstract.}
 Given a symmetrizable generalized Cartan matrix $A$, for any index $k$,
 one can define an automorphism associated with $A,$ of the field
 $\mathbf{Q}(u_1, \cdots, u_n)$ of
 rational functions
 of $n$ independent indeterminates $u_1,\cdots, u_n.$
 It is an isomorphism between two cluster algebras associated to the matrix $A$ (see section 4
 for precise
 meaning).
   When $A$ is of finite type,
  these isomorphisms behave nicely,  they are compatible with the
   BGP-reflection functors
  of cluster categories defined in [Z1, Z2]
  if we identify the indecomposable objects in the categories with cluster
  variables of the corresponding cluster algebras, and they are
  also compatible with
  the "truncated simple reflections" defined in [FZ2, FZ3].
  Using the construction of preprojective or preinjective modules
  of hereditary algebras by Dlab-Ringel [DR] and the Coxeter automorphisms (i.e.,
  a product of these isomorphisms),
  we construct infinitely many cluster variables for cluster algebras
  of infinite type and all cluster variables for finite types.

\medskip


\medskip

}
\end{minipage}
\end{center}
\medskip

\begin{center}

\textbf{1. Introduction}\end{center}

\medskip

Clusters and cluster algebras are defined and studied by Fomin and
Zelevinsky [FZ1-3, BFZ] in order to provide an algebraic framework
for total positivity and canonical bases in semisimple algebraic
groups. Since they appeared, there have been many interesting
connections with other directions [FZ1-3] [BFZ] [CC], amongst them
to representation theory
    of quivers and tilting theory [MRZ] [BMRRT] [BMR] [CC] [CK] [Z1, Z2].

      The connections between representation theory of quivers and
      cluster algebras
      are firstly discovered  by Marsh-Reineke-Zelevinsky through extending
      the well-known Gabriel's Theorem; and then by Buan-Marsh-Reineke-Reiten-Todorov
       who introduced cluster categories, see also [CCS], and related tilting theory with clusters;
       and by some others.
         Since the cluster algebras of finite type are classified by Dynkin
         diagrams [FZ2], there should be some stronger links between
         cluster algebras with
         representation theory of quivers. In [Z1], see also [Z2], we introduced
         the BGP-reflection
         functors in cluster categories by extending the usual
         BGP-functors of module categories. By applying these functors to cluster
         algebras of finite type, we gave a one-to-one correspondence from
         indecomposable objects in cluster
         categories to the almost positive roots of the corresponding
         simple Lie algebras, and then to the set of cluster
         variables of corresponding cluster algebras. This correspondence
         sends basic
         tilting objects
          to clusters.

\medskip

The aim of the paper is to understand clusters and cluster algebras
of any type in terms of representations of quivers based on the
works of [BMRRT], [Z1]. From Dlab-Ringel [DR], any preprojective or
preinjective indecomposable module can be constructed as an image of
projective resp. injective modules under some power of Coxeter
functors. Applying to cluster categories, the indecomposable objects
coming from preprojective or preinjective modules are
 also constructed from the objects $P_k[1]$ by acting some powers of
 Coxeter functors in cluster categories, where $P_k$ is
  an indecomposable projective module and the Coxeter functors in
 cluster categories are defined as a composition of
 BGP-reflection functors introduced in [Z1]. Passing from the set of
 indecomposable
 objects in cluster categories to
 the set of cluster variables, we should get a construction of some
 cluster variables
 from $u_i$ by some automorphisms of cluster algebras induced from Coxeter
 functors in cluster categories. So we should find some
 automorphisms of cluster algebras corresponding to the Coxeter
 functors in the corresponding cluster categories.

  Let $A=(a_{ij})$ be a symmetrizable generalized Cartan matrix
  and $\mathcal{F}=Q(u_1, \cdots, u_n)$ the field of rational functions
   in variables $u_i$ $i=1, \cdots, n.$ For any $k$,
  we define $T_k$ to be the automorphism of $\mathcal{F}$ by
  setting
  $T_k(u_k)=u_k^{-1}(\prod_{a_{ki}<0}u_i^{-a_{ki}}+1),\
  T_k(u_j)=u_j \ \forall j\not= k.$  Let $\underline{u}=(u_1,\cdots, u_n)$ and $B(A)$
  a skew-symmetrizable matrix whose Cartan counterpart is $A$ and such that
  $k$ is a sink in the quiver
   corresponding to $B(A)$. Then $T_k$ is an isomorphism from the cluster algebra associated to
   the initial seed $(\underline{u}, B(A))$ to the cluster algebra associated to
   the initial seed $(\underline{u} ,s_kB(A)).$ When $A$ is of finite type, these automorphisms
   are compatible with the BGP-reflection functors in
  cluster categories when we identify the indecomposable objects
  with cluster variables. By using these isomorphisms, we give a
  construction of cluster variables from the initial cluster. We extend
  the construction to cluster algebras of infinite type.
\medskip

This paper is organized as follows: in Section 2, some basic results
on cluster categories which will be needed later on are recalled. In
Section 3, we recall the BGP-reflection functors in cluster
categories from [Z1, Z2], and extend the result of Dlab-Ringel [DR]
to cluster categories, namely, we prove the indecomposable objects
coming from preprojective or preinjective modules are some powers of
the Coxeter functor on objects $P[1]$ where $P$ is an indecomposable
projective module. In Section 4, for any symmetrizable integer
matrix and any index $i$, we define an automorphism $T_i$ on
$\cal{F}$ and the Coxeter automorphism $T$ as a product of these
$T_i$, and prove that this Coxeter automorphism is a symmetry of the
corresponding cluster algebras. This gives some nice consequences
such as that all cluster variables can be obtained from the initial
cluster by some powers of $T$ when $A$ is of finite type. This
construction is generalized to the infinite types.

 \medskip

 \begin{center}

\textbf{2. Basics on cluster categories.}
\end{center}

\medskip

Let $H$ be a finite-dimensional hereditary algebra over a field
$\mathbf{k},$ with $n$ pairwise non-isomorphic simple modules. Then
there are $n$ pairwise non-isomorphic indecomposable projective
$H-$modules $P_1,\cdots, P_n$. We denote by $\mathcal{D} =
    D^{b}(H)$ the bounded derived category of $H$ with shift
functor $[1]$. For any Krull-Schmidt category $\cal{E}$ [Ri], any
object can be written as a direct sum of indecomposable objects and
such decomposition is unique up to isomorphisms, we will denote by
$\ind\cal{E}$ the full subcategory of representatives of isomorphism
classes of indecomposable objects in $\cal{E}$; depending on the
context we shall also use the same notation to denote the set of
isomorphism classes of indecomposable objects in $\cal{E}$.
\medskip

The cluster category of type $H$ is introduced in [BMRRT], which is
defined to be the factor category $\mathcal{D}/ F$ of $D^b(H)$,
where $F= \tau^{-1}[1]$ and $\tau$ is the Auslander-Reiten
translation in $\cal{D}$. We simply denote the cluster category of
type $H$ by $\mathcal{C}(H)$. This factor category $\mathcal{D}/F$
is a Krull-Schmidt
 triangulated category [K]. The canonical functor $\pi \colon
\mathcal{D} \longrightarrow \mathcal{D}/F: X\mapsto \widetilde{X}$
is a covering functor of triangulated categories, i.e., it sends
triangles to triangles [XZ].  The shift in $\mathcal{D}/F$ is
induced by the shift in $\cal{D}$, and is also denoted by $[1]$. In
both cases we write as usual $\Hom(U,V[1]) = \Ext^1(U,V)$. We then
have
$$\Ext^1_{\mathcal{D}/F}(\widetilde{X},\widetilde{Y}) =
\oplus_{i \in \mathbf{Z}} \Ext^1_{\mathcal{D}}( X, F^iY),$$ where
$X,Y$ are objects in $\cal{D}$ and $\widetilde{X},\widetilde{Y}$ are
the corresponding objects in $\mathcal{D}/F.$

  \medskip

 \textbf{Proposition 2.1.[BMRRT]} Any indecomposable object in $\mathcal{C}(H)$ is
 of the form $\widetilde{M}$ for some indecomposable $H-$module $M$ or
 $\widetilde{P_i[1]}$
 for an indecomposable projective module $P_i, \ 1\le i\le n.$
 \medskip

$H-$mod can be embedded into $D^b(H)$ so that the image of a
 $H-$module is a stalk complex of  degree zero. Passing to the
cluster category $\mathcal{C}(H),$ obj(ind$H)$ can be viewed as a
subset of obj(ind$\mathcal{C}(H))),$ we fix this inclusion in the
rest of paper. Then ind$\mathcal{C}(H)=\mbox{ind}H \cup \{ P_i[1]\ |
\ 1\leq i\leq n \}$ (from now on, the tilde notation for objects in
$\mathcal{C}(H)$ is dropped). For any hereditary algebra $H$, the
indecomposable $H-$modules are either preprojective, or regular, or
preinjective; i.e.,

$$\mbox{ind}H=\mathcal{P}\vee \mathcal{R}\vee\mathcal{I},$$  where
$\mathcal{P}$ (or $\mathcal{I}$) denotes the subcategory of
indecomposable preprojective modules (resp. preinjective modules),
$\mathcal{R}$ denotes the subcategory of indecomposable regular
modules. If $H$ is of finite type, then $\mathcal{R}$ disappears,
and $\mathcal{P}= \mathcal{I}.$ Applying to cluster category, we
have the following:
\medskip

\textbf{Proposition 2.2.} ind$\mathcal{C}(H)=\mathcal{P}\vee\{
P_i[1]\ |\ 1\le i\le n \}
   \vee\mathcal{I}\vee\mathcal{R}.$
 If $H$ is of finite type,  then $\mathcal{R}=\emptyset $,
 $\mathcal{P}= \mathcal{I},$
ind$\mathcal{C}(H)=\mathcal{P}\vee\{ P_i[1]\ |\ 1\le i\le n \}=
   \mathcal{I}\vee\{ P_i[1] \ |\ 1\le i\le n \}.$
\medskip

 \begin{center}

\textbf{3. BGP-reflection functors in cluster categories.}
\end{center}

\medskip

\medskip

Since any hereditary algebra is Morita equivalent to a tensor
algebra of some species of a valued quiver, we will use the
language of valued quivers and their representations. Firstly we
recall some basic notations on representations of valued quivers
from [DR].
\medskip

Let $(\G,\b)$ be a valued graph with $n$ vertices and with an
orientation $\Omega$ (the pair $(\G,\b, \Omega)$ or simply the
pair $(\G, \Omega)$ is called a valued quiver).
   For any  vertex $k\in \G, $ we can define a new orientation
  $s_k\Omega$ of $(\G,\b)$ by reversing the direction of arrows
  along all edges containing $k$. A vertex $k\in \G$ is said to be
  a sink (or a source) with respect to $\Omega$ if there are no arrows
  starting (or ending) at vertex $k$.

  \medskip

  Let $\mathbf{k}$ be a field and $(\G, \Omega)$ a valued quiver.
  From now on, we shall always assume that
  $(\G, \Omega)$ contains no oriented cycles. For any orientation $\Omega$,
  there is an ordering $k_1, \cdots, k_n$ of $\G$ such
that vertex $k_t$ is a sink with respect to the orientation
$s_{k_{t-1}}\cdots s_{k_2}s_{k_1}\Omega $ for all $1\le t\le n.$
This is also equivalent to that
  there is an ordering $k'_1, \cdots, k'_n$ of $\G$ such
that the vertex $k'_t$ is a source with respect to the orientation
$s_{k'_{t-1}}\cdots s_{k'_2}s_{k'_1}\Omega $ for all $1\le t\le n.$
Such orderings are called admissible sequences of sinks or
admissible sequences of sources and an orientation with an
admissible sequence of sinks (equivalently sources) is called an
admissible orientation. It is clear that for any admissible sequence
$k_1, \cdots, k_n$ of sinks or sources, $s_{k_{n}}\cdots
s_{k_2}s_{k_1}\Omega =\Omega.$

 Let
  $\m=(F_i, {}_iM_j)_{i,j\in \G}$ be a reduced $\mathbf{k}-$species of
  type $ \Omega: $ that
  is, for all $i, j \in \G$, $_iM_j$ is an $F_i-F_j-$bimodule,
  where $F_i$ and $F_j$
   are finite extensions of $\mathbf{k}$  and dim$(_{i}M_{j})_{F_j}=d_{ij}$
   and dim$_{\mathbf{k}}F_i=
   \varepsilon_i$. A $\mathbf{k}-$representation $V=(V_i,{}_j\varphi_i)$ of
   $\m$ consists of $F_i-$ vector spaces $V_{i}, i\in \G$, and of an
   $F_j-$linear map $_j\varphi_i: V_i\otimes {}_iM_j\longrightarrow V_j$
   for each arrow $i\longrightarrow j$. Such a representation is called
   finite dimensional if $\sum _{i\in \Gamma }\mbox{dim}_{\mathbf{k}}V_i
   <\infty.$ The
   category of finite-dimensional representations of $\m$ over
   $\mathbf{k}$ is denoted by rep$(\m , \G, \Omega)$.
   If rep$(\m , \G, \Omega)$ contains only finitely many
   indecomposable representations up to isomorphism, then $\G$ is
   called of finite type; otherwise, $\G$ is called of infinite type.
   It was proved by Gabriel [ARS][R] that
   $\G$ is of finite type if and only if $\G$ is a disjoint union of Dynkin
   diagrams.

   \medskip

Now we fix a $\mathbf{k}-$species $\m$ of a given valued quiver
$(\G, \Omega)$. Given a sink, or a source $k$ of the quiver
    $(\G,
\Omega)$,  we recall the reflection functor
    $S^{\pm}_k$:

    $$S^+_k :\  \mbox{rep}(\m ,\G,  \Omega) \longrightarrow  \mbox{rep}(\m
    , \G, s_k\Omega),\  \mbox{    if } \ k  \ \mbox{ is a sink,}$$ or
$$S^-_k :\  \mbox{rep}(\m ,\G,  \Omega) \longrightarrow  \mbox{rep}(\m
,\G, s_k\Omega), \ \mbox{    if } \ k  \ \mbox{ is a source }.$$

\medskip

 We assume $k$ is a sink. For any representation $V=(V_i, \phi _{\alpha})$ of
$(\m , \G, \Omega)$, its image under $S^+_k$ is by definition,
$S^+_kV=(W_i,{}_j\psi_i),$ a representation of $(\m , \G,
s_k\Omega)$, where $W_i=V_i,$ if $i\not=k;$ and $W_k$ is the kernel
in the diagram:
$$\begin{array}{lcccccccc}
(*)&&0&\longrightarrow& W_k&\s{({}_j\chi
_{k})_{j}}{\longrightarrow}& \oplus {}_{j\in \G} V_j\otimes
{}_jM_k &\s{({}_k\phi _{j})_j}{\longrightarrow}&V_k
\end{array}$$
${}_j\psi _i={}_j\phi _{i}$ and ${}_j\psi
 _{k}={}_j\bar{\chi}_{k}: W_k\otimes {}_kM_j\longrightarrow X_j,$
 where ${}_j\bar{\chi}_{k}$ corresponds to ${}_j\chi _{k}$ under
 the isomorphism Hom$_{F_j}(W_{k}\otimes
 {}_kM_j,V_j)\approx\mbox{Hom}_{F_i}(W_k, V_j\otimes
 {}_jM_i ).$

If \textbf{$\alpha$}$=(\alpha _{i}): V\longrightarrow V'$ is a
morphism in $\mbox{rep}(\m ,\G,  \Omega)$, then
$S^+_k$\textbf{$\alpha$}$=\beta=(\beta _i)$, where $\beta _i =\alpha
_i $ for $i\not=k$ and $\beta _k: W_k\longrightarrow W_k'$ is the
restriction of $\oplus _{j\in \G}(\alpha _j\otimes 1)$ given in the
following commutative diagram:

\[ \begin{CD}
0@>>>W_k @>(_j\chi_{k})_{j}>>\oplus_{j\in \G} V_j\otimes {}_jM_k @>(_k\phi_j)_j>>V_k\\
@VVV@VV\beta _k V @VV\oplus _j(\alpha _j\otimes 1)V @VV\alpha _k V \\
0@>>>W_k'@>(_j\chi_{k}')_{j}>>\oplus_{j\in \G} V_j'\otimes
{}_jM_k@>(_k\phi_j')_j>>V_k'
    \end{CD} \]

\medskip

If $k$ is a source, the definition of $S^-_kV$ is dual to that of
$S^+_kV$, we omit it and refer to [DR].

\medskip

 For simplicity, we denote by
 $\cal{H}$ the category rep$(\m, \G, \Omega)$ and by
$\cal{H}'$ the category rep$(\m, \G, s_k\Omega)$, where $k$ is a
sink (or source) of $(\G,\Omega)$. The cluster categories
$D^b(\mathcal{H})/F$, $D^b(\mathcal{H}')/F$ are denoted by
$\mathcal{C}(\Omega)$ and $\mathcal{C}(s_k\Omega)$ respectively.

\medskip

Let $P_i,\  I_i$ (or $P_i'$ , $I'_i$) be the projective, injective
indecomposable representation in $\mathcal{H}$ (resp.
$\mathcal{H}'$) corresponding to the vertex $i\in \G$, and $E_i$
(resp. $E_i'$) the corresponding simple
 representation in $\mathcal{H}$ (resp. $\mathcal{H}'$).
 We denote by $H$ (resp. $H'$) the tensor algebra of $(\m, \G,
\Omega)$ (resp. $(\m, \G, s_k\Omega)$). Note that $H-$mod is Morita
equivalent to $\mathcal{H}$, and if $k$ is a sink (or source), then
$P_k=E_k$ (resp. $I_k=E_k$) is simple projective (resp. injective)
$H-$module.
 \medskip

  Let $T=\oplus _{i\in \G-\{k\}}P_i \oplus \tau ^{-1}P_k$. Suppose
$k$ is a sink, then $T$ is a tilting $H-$module which is called BGP-
or APR-tilting module and $S^+_k=\mbox{Hom}(T,-)$ as functors. The
following theorem was proved in [Z1] (in a more general case).

\medskip

 \textbf{Theorem 3.1.} Let $k$ be a sink (or a source)
  of a valued quiver $(\G, \Omega).$
  Then the BGP-reflection functor $S^+_k$ (resp. $S^-_k$) induces
  a triangle equivalence $R(S^+_k)$
  (resp.,$R(S^-_k)$) from
  $\mathcal{C}(\Omega)$ to
  $\mathcal{C}(s_k\Omega).$  Moreover we have that
  $$\begin{array}{lc}&\\
  &\\
   R(S^+_k)(X)=&\left\{\begin{array}{lccl}P_k'[1],&&&
   X\cong E_k\\
E_k',&&&
   X\cong P_k[1]\\
   P_j'[1],&&&
   X\cong P_j[1], j\not= k \\
S^+_k(X),&&&
   X\in \mbox{ind}H-\{E_k\} \end{array}\right.\\
   &\\&
   \end{array}$$
\medskip

\textbf{Definition 3.2.} The triangle equivalence functor $R(S^+_k)$
from $\mathcal{C}(\Omega)$ to
  $\mathcal{C}(s_k\Omega)$ induced from the reflection functor $S^+_k$
  is called the BGP-reflection functor (for simplicity,
  reflection functor) in $\mathcal{C}(\Omega)$
  at the sink $k$, which is denoted simply by $R^+_k.$ Dually
  for a source $k$, we have the reflection functor $R^-_k$ from
  $\mathcal{C}(\Omega)$ to
  $\mathcal{C}(s_k\Omega)$.

\medskip

Let $k_1, \cdots, k_n$ be an admissible sequence of sinks for the
quiver $(\G,\Omega).$ Set $C^+=R^+_{k_n}\cdots R^+_{k_2}R^+_{k_1},$
the composition of $R^+_{k_i}$. $C^+$ is a self-equivalence of
$\mathcal{C}(\Omega)$, it is called the Coxeter functor in the
cluster category $\mathcal{C}(\Omega).$ For simplicity, we denote
$C^+$ by $C$. The inverse $C^-$ of $C$, which is also called Coxeter
functor in $\mathcal{C}(\Omega)$ is $C^-=R^-_{k_1}R^-_{k_2}\cdots
R^-_{k_n}.$
 If $(\G, \Omega)$ contains no oriented cycles,  then any orientation is admissible,
  and there are exactly two orientations on $\G$ such that any
 vertex is sink or source. In this case we use $\Omega_0$ and $\Omega_0'$ to
  denote these two distinct orientations on $\G.$
\medskip

\textbf{Theorem 3.3.} If $\G$ is of finite type, then
 $\mbox{ind}\mathcal{C}(\Omega)=\mathcal{P}\vee\{ P_i[1] \ |\ i\in \G \}
   =\mathcal{I}\vee\{ P_i[1] \ |\ i\in \G \};$ if $\G$ is of
   infinite type, then
$\mbox{ind}\mathcal{C}(\Omega)=\mathcal{P}\vee\{ P_i[1] \ |\ i\in
\G \}
   \vee\mathcal{I}\vee\mathcal{R};$
 Moreover, we have that $$\begin{array}{ll}\mathcal{P}=\{ C^{-m}P_i[1]\ |
   \ \mbox{} m> 0, i\in \G\ \}, &  \mathcal{I}=\{ C^{m}P_i[1]\ |
   \ \mbox{ } m> 0, i\in \G\ \}.\end{array}$$

\textbf{Proof.} The first part follows from Proposition 2.2. We
prove the second part. For any indecomposable projective
representation $P$ in $\mathcal{H}$, we have that
$\tau^{-1}P[1]\simeq P$ in $\mathcal{C}(\Omega)$. It follows that
$C^{-1}P[1]\simeq P$ in $\mathcal{C}(\Omega)$. From [DR], we know
that for any preprojective indecomposable module $M$, there are an
indecomposable projective module $P$ and an integer $m\ge 0$, such
that $M=C^{-m}P$. Therefore we have that $M=C^{-(m+1)}P[1].$ This
finishes the proof of the description of $\mathcal{P}$. Dually, for
any preinjective $H-$module $N$, there are an indecomposable
injective module $I$ and an integer $m\ge 0$ such that $N=C^{m}I$.
We also have that $I_i= \tau P_i[1]$ in $\mathcal{C}(\Omega)$ for
$i\in \G$ (since in derived category $D^b(\mathcal{H})$, $\tau
P_i[1]=I_i$). Therefore we have that $N=C^{m+1}I$. This finishes the
proof for $\mathcal{I}.$ The proof is finished.
\medskip

For this reason, we denote the union $\mathcal{P}\vee\{ P_i[1] \}
   \vee\mathcal{I}$ by $\mathcal{PI}(\Omega).$  Note that
   $$\mathcal{PI}(\Omega)=\{ C^{m}P_k[1]\ \ |
   \  m\in \mathbf{Z}, 1\le k\le n\ \}.$$

   If $\G$ is of finite type, then ind$\mathcal{C}(\Omega)
   =\mathcal{PI}(\Omega),$ otherwise ind$\mathcal{C}(\Omega)
   =\mathcal{PI}(\Omega)\vee \mathcal{R}.$

\medskip

\textbf{Corollary 3.4.} If $\G$ is of finite type, then for any
orientation $\Omega$ on $\G$, ind$\mathcal{C}(\Omega)=\{
C^m(P_k[1])\ |\ m\in \mathbf{Z}, \ 1\le k\le n\ \}=\{ C^m(P_k[1])\
|\ m\in \mathbf{Z}, \ m\ge 0, 1\le k\le n\ \}=\{ C^{-m}(P_k[1]) \ |\
m\in \mathbf{Z}\
 m\ge 0 , 1\le k\le n\ \}.$
\medskip

\textbf{Example 3.5.} Let $\G$ be $B_2: 2\s{(1,2)}{-}1. $  We give
it an
 orientation $\Omega:  2 \s{(1,2)}{\longrightarrow} 1.$ The AR-quiver of the
cluster category $\mathcal{C}(\Omega)$ is the following:

 \setlength{\unitlength}{1cm}
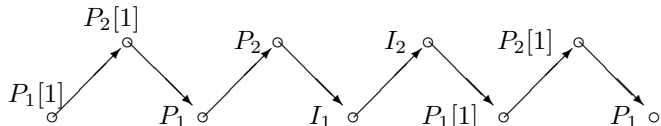
\begin{picture}(4,2)
\mu(0.1,0.1)(2,0){4}{\vector(1,1){0.9}}

\mu(1.1,1.1)(2,0){4}{\vector(1,-1){0.85}}
 \mu(0,0)(2,0){5}{$\c$}
\mu(1,1)(2,0){4}{$\c$}

\put(-0.5,0.3){${\small P_1[1]}$} \put(0.5,1.3){$P_2[1]$}
 \put(1.5,0){$P_1$} \put(2.5,1){$P_2$}

\put(3.5,0){$I_1$} \put(4.5,1){$I_2$} \put(5,0){$P_1[1]$}
 \put(6,1){$P_2[1]$}
\put(7.5,0){$P_1$}
\end{picture}

with the valuation $(1,2)$ on all arrows like $\nearrow$, and the
valuation $(2,1)$ on all arrows like $\searrow$.

\medskip

\medskip

The Coxeter functor $C=R_2^+R_1^+,$ and we have
$P_1=C^{-1}P_1[1],\ P_2=C^{-1}P_2[1],\ I_1=C^{-2}P_1[1] ,\
I_2=C^{-2}P_2[1]; $  also $P_1=C^{2}P_1[1],\ P_2=C^{2}P_2[1],\\
I_1=CP_1[1],\  I_2=CP_2[1]. $
\medskip

 For a valued graph $\G$, we denote by $\Phi$ the set of roots of the
corresponding Kac-Moody Lie algebra. Let $\Phi_{\geq -1}$ denote
the set of almost positive roots, i.e.\ the positive roots
together with the negatives of the simple roots.
  Let $s_i$ be the Coxeter generator of the
Weyl group of $\Phi$ corresponding to $i\in \G $. We recall from
[FZ3] that the "truncated reflections" $ \sigma_i$ of $\Phi_{\geq
-1}$ are defined as follows:
$$  \sigma_i(\alpha)=\left\{ \begin{array}{ll} \alpha & \alpha=-\alpha_j,\ j\not=i \\ s_i(\alpha) & \mbox{otherwise.} \end{array}\right.$$

\medskip

  On the one hand, when $\G$ is of finite type, there is a bilinear form $(-\parallel-)$ on $\Phi_{\geq
-1}$ which is called the "compatibility degree" of $\Phi_{\geq -1}$
(for details, we refer to [FZ3, FZ2]). $\alpha, \ \beta \in
\Phi_{\geq -1}$ are called compatible if $(\alpha ||\beta)=0$. Any
maximal mutually compatible subset is called a cluster of
$\Phi_{\geq -1}$. It was proved in [FZ3] that any cluster in
$\Phi_{\geq -1}$ contains $n$ elements, where $n$ is the number of
simple roots of $\Phi_{\geq -1}$.

On the other hand, in the cluster category $\mathcal{C}(H)$, there
is a tilting machinery. An object $T$ is called tilting if
Ext$_{\mathcal{C}(H)}(T,T)=0$ and it has a maximal number of
non-isomorphic indecomposable direct summands. A multiplicity-free
tilting object is called a basic tilting object. In the following we
assume tilting objects are always basic. Any tilting object contains
$n$ indecomposable direct summands.

We have seen that ind$H\subset \mbox{ind}\mathcal{C}(H),$ and
$\Phi\subset\Phi_{\ge-1}.$ The well-known Gabriel's Theorem gives a
one-to-one correspondence from ind$H$ of hereditary algebra of
finite type to the root system $\Phi^+$ of the corresponding simple
Lie algebra by taking dimension vectors of modules. This
correspondence was generalized to cluster categories of finite type,
which induces a bijection between the set of tilting objects to the
set of clusters in $\Phi _{\geq -1},$ in the simple-laced case in
[BMRRT], and to all Dynkin cases in [Z1] (see Proposition 3.7.
latter for precise meaning). In fact this map can be defined for any
cluster category (finite and infinite types) as follows: for any
$X\in \mbox{ind}(\mbox{mod}H\vee H[1]),$
$$\gamma_{\Omega}(X)=\left\{
\begin{array}{lrl}\mathbf{dim}X & \mbox{ if } & X\in \mbox{ind}H;\\
&&\\
-\mathbf{dim}E_i& \mbox{ if } &X=P_i[1],\end{array}\right.$$ where
$\mathbf{dim}X$ denotes the dimension vector of the representation
$X$. In general, this map $\gamma_{\Omega}:\
\ind\mathcal{C}(\Omega)\rightarrow
 \Phi_{\geq -1} $ is not injective, but it is surjective in
all cases, and is a bijection in the finite type case.

\medskip

Let $\Phi_{\geq -1} '$ be the subset of $\Phi_{\geq -1} $ consisting of
 $\gamma_{\Omega}(X)$ for all $X\in \mathcal{PI}(\Omega)$. Then the restriction of $\gamma _{\Omega}$ to
 $\mathcal{PI}(\Omega)$ is a bijection from  $\mathcal{PI}(\Omega)$ to $\Phi_{\geq -1} '$
 by [DR], [Kac], this map is also denoted by $\gamma _{\Omega}$.

 \medskip

 When $\G$ is of
finite type, one can choose a skew-symmetrizable integer matrix
$B$ from $\G$ such that the cluster variables of type $\G$ are in
one-to-one
 correspondence with the elements of $\Phi_{\geq -1}$ (compare
[FZ2]). For any orientation $\Omega$, ind$\mathcal{C}(\Omega)$ is
in one-to-one correspondence
 with $\Phi_{\geq -1}$ [BMRRT] [Z1]. We will relate these two results
 and generalize partially
 these one-to-one correspondences
 to infinite type in the next section.

\medskip

 By using Theorem 3.1, one gets the following commutative diagram which
 explains that $R_{k}^{\pm}$ is the realization of the
 "truncated reflection" $\sigma _k$ (for proof, we refer to [Z1, Z2]).

 \medskip

 \textbf{Proposition 3.6.}  Let $k$
be a sink (or a source) of a valued quiver $(\G, \Omega)$. Then we
have the commutative diagram:
\[ \begin{CD}
\mbox{ind}\mathcal{C}(\Omega) @>R_{k}^+
>(R_{k}^-)> \mbox{ind}\mathcal{C}(s_k\Omega)
\\
@V\gamma _{\Omega} VV  @VV\gamma  _{s_k\Omega}V  \\
\Phi_{\geq -1} @>\sigma _k>> \Phi_{\geq -1}
\end{CD} \].

\medskip

The following result is proved  for simply-laced Dynkin diagram in [BMRRT], and is generalized to all Dynkin diagram in [Z1].
\medskip

 \textbf{Proposition 3.7 [BMRRT] [Z1].} Let $(\G, \Omega)$ be any valued
 Dynkin quiver. Then the one-to-one correspondence $\gamma
 _{\Omega}$ sends tilting objects of $\mathcal{C}(\Omega)$
 to clusters in $\Phi_{\geq -1}.$

\newpage

 \begin{center}

\textbf{4. Coxeter automorphisms of cluster algebras.}
\end{center}

 \medskip

We recall some basic notation on cluster algebras which can be found
in the series of papers by Fomin and Zelevinsky [FZ1,FZ2, FZ3, BFZ].
The cluster algebras we deal with in this paper are defined on a
trivial semigroup of coefficients, since it is enough for the
connection with representation theory of quivers [BMRRT]. These
cluster algebras are called reduced cluster algebras in [CC]. We
will call these algebras just cluster algebras.
\medskip

The definition is as follows: Let
$\mathcal{F}=\mathbf{Q}(u_1,u_2,\cdots, u_n)$ be the field of
rational functions in indeterminates  $u_1,u_2,\cdots, u_n.$
 Set $\underline{u}=(u_1,u_2,\cdots,
u_n).$ Let $B=(b_{ij})$ be an $n\times n$ skew-symmetrizable integer
matrix. A pair $(\underline{x}, B)$, where
$\underline{x}=(x_1,x_2,\cdots, x_n)$ is a transcendence base of
$\mathcal{F}$ and where $B$ is an $n\times n$ skew-symmetrizable
integer matrix, is called a seed. Fix a seed $(\underline{x}, B)$
and an element $z$ in the base $\underline{x}$. Let $z'$ in
 $\mathcal{F}$ be such that
 $$zz'=\prod_{b_{xz}>0}x^{b_{xz}}+\prod_{b_{xz}<0}x^{-b_{xz}}.$$

 Now, set $\underline{x}':=\underline{x}-\{z\}\bigcup \{z'\}$
  and $B'=(b'_{xy})$ such that
  $$b'_{xy}=\left\{\begin{array}{lccccl}-b_{xy}&&&&&\mbox{if } x=z
  \mbox{ or }y=z,\\
  b_{xy}+1/2(|b_{xz}|b_{zy}+b_{xz}|b_{zy}|)&&&&&
  \mbox{otherwise.}\end{array}\right.$$

  The pair $(\underline{x}', B')$ is called the mutation
  of the seed $(\underline{x}, B)$ in direction $z$, it is also a seed.
  The "mutation equivalence $\approx $" is an equivalence relation
  on the set of all seeds generated
  by $(\underline{x}, B)\approx(\underline{x}', B')$  if
  $(\underline{x}', B')$ is a mutation of $(\underline{x}, B)$.

The cluster algebra $\mathcal{A}(B)$ associated to the
skew-symmetrizable matrix $B$ is by definition the subalgebra of
$\mathcal{F}$ generated by all $\underline{x}$ such that
$(\underline{x}, B')\approx(\underline{u}, B).$ Such
$\underline{x}=(x_1,x_2,\cdots, x_n)$ is called a cluster of the
 cluster algebra $\mathcal{A}(B)$ or simply of $B$, and any
 $x_i$ is called a cluster variable. If the set $\chi$ of all
 cluster variables is finite, then the cluster
  algebra $\mathcal{A}(B)$ is said to be of finite type.
\medskip

Let $\underline{x}$ be a seed. The Laurent phenomenon, see [FZ1],
asserts that any cluster variables are Laurent polynomials with
integer coefficients in variables $x_1, \cdots, x_n.$ It implies
that $\mathcal{A}(B)\subset \mathbf{Z}[x_1^{\pm}, \cdots,
x_n^{\pm}].$

Fix any integer square matrix $B=(b_{ij})$. Its Cartan counterpart
is by definition, a generalized Cartan matrix $A=A(B)=(a_{ij})$ of
the same size defined by
$$a_{ij}=\left\{\begin{array}{lccl}2&&&
\mbox{if }i=j\\
-|b_{ij}|&&& \mbox{if } i\not=j.\end{array}\right.$$
\medskip

\textbf{Theorem 4.1.}[FZ2]  A cluster algebra $\mathcal{A}$
 is of finite type if and only if there is a seed $(\underline{x}, B)$ of
 $\mathcal{A}$ such that the Cartan counterpart of the matrix $B$ is a
 Cartan matrix of finite type.

 \medskip

 For any $\alpha \in \Phi$, we write $\alpha$ as a sum of simple roots
 $\alpha=\sum_{i\in I}a_i\alpha _i$, then we use $u^{\alpha}$ to denote $\prod_{i\in
 I}u_i^{a_i}.$
 \medskip

\textbf{Theorem 4.2.}[FZ2]  Fix a Dynkin diagram $\G$ and a distinguished seed $(\underline{u},B)$. Then there
exists a bijection $$P: \Phi_{\ge-1} \longrightarrow \chi_{\G}: \
\alpha \mapsto u[\alpha]=\frac{P_{\alpha}(u)}{u^{\alpha}},$$ where
$P_{\alpha}(u)$ is a polynomial with nonzero constant term. Under
this
 correspondence, $-\alpha _i$ corresponds to $u_i$ and clusters in
$\Phi_{\ge-1}$ correspond to clusters of the corresponding
 cluster algebra $\mathcal{A}.$

\medskip

\medskip

\textbf{Corollary 4.3.}  Let $\G$ be a Dynkin graph with the
orientation $\Omega_0$ (i.e., such that any vertex is a sink or
source). Then
 the composition (denoted by $\phi_{\Omega _0}$)
 of $\gamma _{\Omega _0}$ and $P$ gives a one-to-one
 correspondence between ind$\mathcal{C}(\Omega_0)$ and $\chi_{\G}$.
  Under this correspondence, $P_i[1]$ corresponds to $u_i$, and
  tilting objects correspond to clusters.

\medskip

\textbf{Proof.} This is a consequence of Proposition 3.7. and
Theorem 4.2.

\medskip

Note that when $\G$ is a simply-laced Dynkin diagram,  Corollary
4.3. is also proved in [CC] (compare Theorem 3.4. there) in which
the correspondence is given in explicit expressions of
 indecomposable objects by Laurent polynomials of $u_1,\ \cdots,
u_n.$ In the following, for any orientation, we will give an
explicit
 one-to-one correspondence from cluster categories to the set of
cluster variables in a different spirit, which works in simply-laced
case and non-simply-laced case. Firstly we need to define some
isomorphisms between cluster algebras.

 Given a generalized Cartan matrix $A$ of size $n\times n,$ its
Coxeter graph $\Delta$ is by definition, a valued graph consisting
of $n$ vertices, named $1,\ 2, \cdots, \ n$, and edges $i- j$ with a
valuation $(a_{ij}, a_{ji})$ if $a_{ij}\not=0.$ If $A$ is a Cartan
matrix of finite type, then its Coxeter graph is a tree.
\medskip

  Let $B=(b_{ij})$
be a skew-symmetrizable integer matrix and $A=A(B)$, its Cartan
counterpart. Then we say that $B$ and $A$ form a matched pair
$(B,A).$ Note that $A$ is a generalized Cartan matrix and there are
different skew-symmetrizable
 matrices $B$ and $B'$ with the same Cartan counterpart $A.$

 For any matched pair $(B,A)$ of matrices with $B$ a skew-symmetrizable
 integer matrix, one can
give an orientation $\Omega $ on its Coxeter graph $\Delta$ as
follows: if $b_{ij}>0,$ then there is an arrow $i\longrightarrow j$.
With this orientation, $(\Delta, \Omega)$ becomes a valued quiver.
This quiver is called the quiver of $(B,A)$. Note that if we can
choose $B=B(A)$
 with the property $"b_{ij}b_{ik}\ge 0,\ $ for all $i, j, k,"$ then
 the quiver of $(B,A)$ is such that any vertex is
 a sink or a source. If a Coxeter graph is a tree, then we can choose
 such an orientation, this orientation was considered in [FZ2]. More
 generally, the orientation on a Coxeter graph corresponds to the
 skew-symmetrizable  matrix $B=(b_{ij})$ such
that $A=A(B)$ is as indicated in the next lemma.

 \medskip

\textbf{Lemma 4.4.} Fix a Coxeter diagram $\Delta$, equivalently, a
generalized Cartan matrix $A$. Then the orientations of $\Delta$ are
in bijection with
 the matched pairs $(B, A)$. Moreover, the orientation contains
 an orientated cycle if and only if there is a sequence of indices
 $i_1, i_2, \cdots i_t$ such that $b_{i_1, i_2}, \ b_{i_2, i_3},\
 \cdots , b_{i_{t-1},i_t},\ b_{i_t,i_1}$ are positive integers.
 \medskip

 \textbf{Proof.} For any matched pair $(B,A)$, we can define an
 orientation of $\Delta$ as above. Conversely, for any orientation
  of $\Delta$, the matrix $B=(b_{ij})$ can be defined uniquely in
 the way:
$$b_{ij}=\left\{ \begin{array}{lccl}0&&&
\mbox{if }i=j\\
|a_{ij}|&&& \mbox{if } i\rightarrow j\\
a_{ij}&&& \mbox{if } j\rightarrow i. \end{array}\right.$$ The final
statement follows from the correspondence between orientations and
$B$. The proof is finished.
\medskip

\textbf{Definition 4.5.} Let $A$ be a generalized Cartan matrix and
$B$ one of the skew-symmetrizable  matrix with $A=A(B).$  For any
index $i$, we define an automorphism $T_i$ of
$\mathcal{F}=\mathbf{Q}(u_1,\cdots,u_n)$ by defining the images of
the indeterminates $u_1, \cdots, u_n$ as follows:

$$T_i(u_j)=\left\{ \begin{array}{lccccl}u_j&&&&&\mbox{if } j\not=i,\\
&&&&&\\
  \frac{\prod_{a_{ik}<0}u_k^{-a_{ik}}+1}{u_i}&&&&&
 \mbox{if } j=i.\end{array}\right.$$

  \medskip

  It is easy to check that all $T_i$ are involutions of $\mathcal{F}$, i.e.,
  $T_i^2=\mbox{id}_{\mathcal{F}}.$

\medskip

From the definition, all $T_i$ are independent of $B$, and only
depend on the matrix $A$.
\medskip

In the rest of paper, we study the properties of automorphisms $T_i$
with respect to clusters and cluster algebras. Since any orientation
 of $\G$ is admissible, for such an orientation $\Omega$ on $\G$, by Lemma
4.4, we have a pair $(B,A)$ with $A$ a generalized Cartan matrix,
the
  matrix $B$ corresponding to $\Omega$ is sometimes denoted by
  $B_{\Omega}.$ The
cluster algebra $\mathcal{A}$ associated with the seed
$(\underline{u}, B_{\Omega})$ is called the cluster algebra
associated with the orientation $\Omega $, and is denoted by
$\mathcal{A}_{\Omega};$
 the set of cluster variables of $\mathcal{A}_{\Omega}$ is denoted by
$\chi_{\Omega}.$ By [FZ2], $\mathcal{A}_{\Omega}$ is isomorphic to
$\mathcal{A}_{\Omega _0}.$

\medskip

\textbf{Remark 4.6.} $\chi_{\Omega}$ is different from
$\chi_{\Omega_0}$ in general, see Example 4.9. below.

\medskip

  In the following, we prove that for any orientation of a
  Dynkin graph $\G$ there is a bijection $\phi _{\Omega}$ from
ind$\mathcal{C}(\Omega)$ to $\chi_{\Omega}$ such that $P_i[1]$
corresponds to $u_i$, which generalizes Corollary 4.3.
 \medskip

\medskip

\textbf{Theorem 4.7.} Let $k$ be a sink (or source) in $\Omega$ on
 a Dynkin diagram $\G.$ Then (1). there is a bijection $\phi _{\Omega}$ from
ind$\mathcal{C}(\Omega)$ to $\chi_{\Omega}$ such that $P_i[1]$
corresponds to $u_i$, for any $i\in \G;$ inducing a one-to-one correspondence
between basic tilting objects and clusters.

(2). $T_k$ sends cluster variables and
 clusters in $\chi_{\Omega}$ to those in
$\chi_{s_k\Omega}.$

(3). $T_k$ is induced from the reflection functor $R^+_k$
indicated in the following commutative diagram:
\[ \begin{CD}
\mbox{ind}\mathcal{C}(\Omega) @>R^+_k
>> \mbox{ind}\mathcal{C}(s_k\Omega)\\
@V\phi_{\Omega} VV  @VV\phi _{s_k\Omega}V  \\
\chi_{\Omega} @>T_k>> \chi_{s_k\Omega}
\end{CD}. \]

(4). $T_k$ induces an isomorphism from the cluster algebra
$\mathcal{A}_{\Omega}$ to the cluster algebra
 $\mathcal{A}_{s_k\Omega}$ ($T_k$ induces a so-called strongly isomorphism from
$\mathcal{A}_{\Omega}$ to
 $\mathcal{A}_{s_k\Omega}$).
\medskip
\medskip

 \textbf{Proof.} By definition, $\chi_{\Omega}$
  and $\chi_{s_k\Omega}$ are the sets of cluster variables of
  the initial seeds $(\underline{u}, B_{\Omega})$ and
  $(\underline{u}, B_{s_k\Omega})$ respectively. Let
  $\underline{u}'_k=  (u_1, \cdots, u_{k-1},\\
   u'_k, u_{k+1}, \cdots, u_n)$
  be the cluster of $(\underline{u}, B_{s_k\Omega})$ obtained
  by  mutation once in
  direction $k.$ Then $u_k'u_k=\prod_{a_{ik}<0}u_i^{-a_{ik}}+1$ since
   $k$ is a source of $s_k\Omega$, hence
   $u_k'=\frac{\prod_{a_{ik}<0}u_i^{-a_{ik}}+1}{u_k},$ and the matrix
  $B_{s_k\Omega}'$ after this mutation is $B_{\Omega}$. Therefore
  we have that $T_k(\underline{u})=\underline{u}'_k$ and
  the automorphism $T_k$ of $\mathcal{F}$
  sends the seed $(\underline{u}, B_{\Omega})$ to the seed
  $(\underline{u}'_k, B_{\Omega}).$  Dually $T_k$ sends the seed
  $(\underline{u}'_k, B_{\Omega})$ to the seed
  $(\underline{u}, B_{\Omega}).$ Therefore, $T_k$ sends
  $\chi _{\Omega}$ to $\chi _{(\underline{u}'_k, B_{\Omega})}$, the latter is the
  set of cluster variables associated to the initial seed $(\underline{u}'_k, B_{\Omega})$.
   Since the seed $(\underline{u}'_k, B_{\Omega})$ is obtained by
   mutation from $(\underline{u}, B_{s_k\Omega}),$ i.e.,
  $(\underline{u}'_k, B_{\Omega})
  \s{mutation}{\approx}(\underline{u}, B_{s_k\Omega}),$
  $\chi _{(\underline{u}'_k, B_{\Omega})}=
  \chi_{s_k\Omega}. $  Hence $T_k$ induces a
  bijection from $\chi_{\Omega}$
  to $\chi_{s_k\Omega}$.  For any seed $(\underline{x}, B)$ which is
  mutation equivalent to the initial seed $(\underline{u},
  B_{\Omega})$, denote by $T_k(\underline{x})$ the vector $(T_k(x_1), \cdots, T_k(x_n))$,
  then
  $(T_k(\underline{x}), B)$ is a seed which is mutation equivalent to the seed
  $(\underline{u}'_k, B_{\Omega})$, and then it is mutation equivalent to the initial seed
  $(\underline{u}, B_{s_k\Omega}).$ This implies that $T_k(\underline{x})$ is a cluster of
  $(\underline{u}, B_{s_k\Omega}).$ Then $T_k$ sends clusters in $\chi_{\Omega}$ to clusters
  in $\chi_{s_k\Omega}.$
   This proves part (2).
  Since $\Omega$ and $\Omega_0$ are
  orientations on $\G$, there is a sequence of vertices $i_1, \cdots, i_t$ such that
  $s_{i_t}\cdots s_{i_1}\Omega_0=\Omega$, where $i_j$ is a sink of
   $s_{i_{j-1}}\cdots s_{i_1}\Omega _0$ for any $j$.  It follows
   from part 2 proved above that
  $T_{i_t}\cdots T_{i_1}\chi_{\Omega_0}=\chi_{\Omega}.$ Combining with Corollary
  4.3, we have the bijection $\phi_{\Omega}$ (obtained
  by taking $\phi_{\Omega}=
 T_{i_t}\cdots T_{i_1} \phi_{\Omega_0}(R^+_{i_t}\cdots R^+_{i_1})^{-1}$) from
  ind$\mathcal{C}(\Omega)$ to $\chi_{\Omega}$ which satisfies
  the commutative diagram:

\[ \begin{CD}
\mbox{ind}\mathcal{C}(\Omega_0) @>R^+_{i_t}\cdots R^+_{i_1}
>> \mbox{ind}\mathcal{C}(\Omega)\\
@V\phi_{\Omega_0} VV  @VV\phi _{\Omega }V  \\
\chi_{\Omega_0} @>T_{i_t}\cdots T_{i_1}>> \chi_{\Omega}
\end{CD} \]

Since all isomorphisms $T_i$ send clusters to clusters, all $R^+_i$
are triangle equivalences sending basic tilting objects to basic
tilting objects, and $\phi _{\Omega_0}$ send basic tilting objects
to clusters (by Corollary 4.3.), we have that $\phi _{\Omega}$ sends
basic tilting objects to clusters. This proves part (1).

Let $i_1, \cdots, i_t$ be the vertices of $\G$ such that
$s_{i_t}\cdots s_{i_1}\Omega_0=\Omega$ with $i_j$ a sink of
   $s_{i_{j-1}}\cdots s_{i_1}\Omega _0$ for any $j$.  By part (1), we have that
   $\phi_{\Omega}=
 T_{i_t}\cdots T_{i_1} \phi_{\Omega_0}(R^+_{i_t}\cdots R^+_{i_1})^{-1}$ and
  $\phi_{s_k\Omega}= T_k T_{i_t}\cdots T_{i_1}
  \phi_{\Omega_0} (R^+_kR^+_{i_t}\cdots R^+_{i_1})^{-1}.$
  Then $T_k\phi_{\Omega}=T_k T_{i_t}\cdots T_{i_1} \phi_{\Omega_0}
  \newline (R^+_{i_t}\cdots R^+_{i_1})^{-1}=T_k T_{i_t}\cdots T_{i_1} \phi_{\Omega_0}
  (R^+_{i_1}\cdots R^+_{i_t}R^+_k)R^+_k
  =\phi_{s_k\Omega}R^+_k.$ This proves (3).

  Since cluster algebras $\mathcal{A}_{\Omega}$ and $\mathcal{A}_{s_k\Omega}$ are generated
   as subalgebras of $\mathcal{F}$ by $\chi_{\Omega}$  and
   $\chi_{s_k\Omega}$ respectively, $T_k$ induces an isomorphism from
   $\mathcal{A}_{\Omega}$ to $\mathcal{A}_{s_k\Omega}$. This proves
   the final part. The whole proof is completed.

\medskip

\textbf{Remark 4.8.} By using the commutative diagram in Theorem 4.7., one can get that the bijection $\phi _{\Omega}$ is of the
form: $C^k(P_i[1])\mapsto T^k(u_i) $ for any $k\in\mathbf{Z}, \ i=1,\cdots, n.$

\medskip

\textbf{Remark 4.9.} In general, for any diagram $\G$ (finite type
or infinite type), $T_k$ sends cluster variables and
 clusters in $\chi_{\Omega}$ to those in
$\chi_{s_k\Omega}$ and $T_k$ induces an isomorphism from
$\mathcal{A}_{\Omega}$ to $\mathcal{A}_{s_k\Omega}$. The proof for
this general result is same as that for part (2) (4) of Theorem 4.7.

     \medskip

\textbf{Example 4.10.} Let $\G$ be $A_3: 3-2-1$ with Cartan
     matrix $A=\left(\begin{array}{ccc}2&-1&0\\
     -1&2&-1\\
     0&-1&2\end{array}\right)$. Let $\Omega$ be an orientation of
     $\G:\ 3\longrightarrow2\longrightarrow1.$
      Then the corresponding skew-symmetrizable matrix $B$ is
      $\left(\begin{array}{ccc}0&-1&0\\
     1&0&-1\\
     0&1&0\end{array}\right).$  The AR-quiver of $\mathcal{C}(\Omega)$ has
     the following shape:

 \setlength{\unitlength}{1cm}
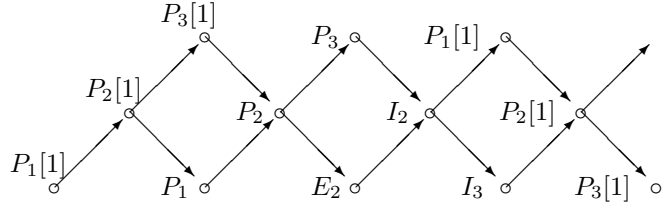
\begin{picture}(4,3)
\mu(0.1,0.1)(2,0){4}{\vector(1,1){0.9}}
\mu(1.1,1.1)(2,0){4}{\vector(1,1){0.9}}
\mu(1.1,1.1)(2,0){4}{\vector(1,-1){0.85}}
\mu(2.1,2.1)(2,0){3}{\vector(1,-1){0.85}} \mu(0,0)(2,0){5}{$\c$}
\mu(1,1)(2,0){4}{$\c$} \mu(2,2)(2,0){3}{$\c$}

\put(-0.5,0.3){${\small P_1[1]}$} \put(0.5,1.3){$P_2[1]$}
\put(1.5,2.3){$P_3[1]$} \put(1.5,0){$P_1$} \put(2.5,1){$P_2$}
\put(3.5,2){$P_3$}

\put(3.5,0){$E_2$} \put(4.5,1){$I_2$} \put(5.5,0){$I_3$}
\put(5,2){${\small P_1[1]}$} \put(6,1){$P_2[1]$}
\put(7,0){$P_3[1]$}
\end{picture}
\medskip

The correspondence $\phi _{\Omega}$ from ind$\mathcal{C}(\Omega)$ to
the set $\chi_{\Omega}$ is indicated as follows:
$$\begin{array}{l}P_i[1]\mapsto u_i,\mbox{ for } i=1,\ 2,\ 3; \
P_1\mapsto \frac{1+u_2}{u_1}; \ P_2\mapsto
\frac{u_1+u_3+u_2u_3}{u_1u_2}; \\
P_3\mapsto \frac{u_1+u_3+u_1u_2+u_2u_3}{u_1u_2u_3}; \
  E_2\mapsto
\frac{u_1+u_3}{u_2};\  I_2\mapsto
   \frac{u_1+u_3+u_1u_2}{u_2u_3}; \
I_3\mapsto \frac{1+u_1}{u_3}.\end{array}$$ The cluster algebra
$\mathcal{A}_{\Omega}$ is the $\mathbf{Q}-$subalgebra of
$\mathcal{F}$ generated by all cluster variables above.

\medskip

If we reflect the orientation $\Omega$ at vertex $1$, we get the
quiver $(\G, s_1\Omega):\ \ 3\longrightarrow2\longleftarrow 1.$ It
corresponds to the skew-symmetrizable matrix
 $B=\left(\begin{array}{ccc}0&1&0\\
     -1&0&-1\\
     0&1&0\end{array}\right).$  The AR-quiver of
     $\mathcal{C}(s_1\Omega)$ has
     the following shape:

 \setlength{\unitlength}{1cm}
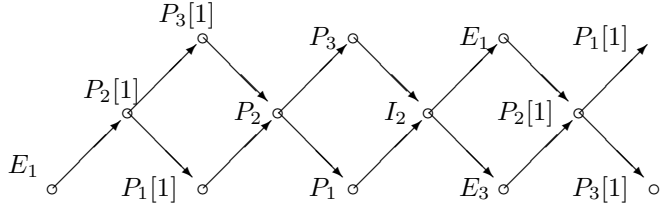
\begin{picture}(4,3)
\mu(0.1,0.1)(2,0){4}{\vector(1,1){0.9}}
\mu(1.1,1.1)(2,0){4}{\vector(1,1){0.9}}
\mu(1.1,1.1)(2,0){4}{\vector(1,-1){0.85}}
\mu(2.1,2.1)(2,0){3}{\vector(1,-1){0.85}}
 \mu(0,0)(2,0){5}{$\c$}
\mu(1,1)(2,0){4}{$\c$} \mu(2,2)(2,0){3}{$\c$}

\put(-0.5,0.3){$E_1 $} \put(0.5,1.3){$P_2[1]$}
\put(1.5,2.3){$P_3[1]$} \put(1,0){$P_1[1]$}

\put(2.5,1){$P_2$} \put(3.5,2){$P_3$}

\put(3.5,0){$P_1$} \put(4.5,1){$I_2$} \put(5.5,0){$E_3$}
\put(5.5,2){$ E_1$} \put(6,1){$P_2[1]$} \put(7,0){$P_3[1]$}
\put(7,2){$P_1[1]$}
\end{picture}

\medskip

In this case, the correspondence $\phi _{s_1\Omega}$ from
ind$\mathcal{C}(s_1\Omega)$ to the set $\chi_{s_1\Omega}$ is as
 follows: $$\begin{array}{l} P_i[1]\mapsto u_i, \mbox{ for } i=1,\
2,\ 3;  \ P_1\mapsto \frac{1+u_2+u_1u_3}{u_1u_2};\ P_2\mapsto
\frac{1+u_1u_3}{u_2}; \\ P_3\mapsto \frac{1+u_2+u_1u_3}{u_2u_3};\
E_1\mapsto \frac{1+u_2}{u_1};\ I_2\mapsto
\frac{1+2u_2+u_2^2+u_1u_3}{u_1u_2u_3};\ E_3\mapsto
\frac{1+u_2}{u_3}\end{array}.$$

The corresponding cluster algebra $\mathcal{A}_{s_1\Omega}$ is the
$\mathbf{Q}-$subalgebra of $\mathcal{F}$ generated by all cluster
variables above. $\chi_{\Omega}\not=\chi_{s_1\Omega}.$ It is easy to
see that $\phi_{s_1\Omega}(R^+_1X)=T_1(\phi_{\Omega}(X))$ for all
$X\in \mbox{ind}\mathcal{C}(\Omega)$. Therefore the isomorphism
$T_1$ induces an isomorphism from the cluster algebra
$\mathcal{A}_{\Omega}$ to the cluster algebras $
\mathcal{A}_{s_k\Omega}.$

\medskip

Let $\Omega$ be an orientation of $\G$ and $C=R^+_{k_n}\cdots
R^+_{k_2}R^+_{k_1}$ the corresponding Coxeter
 functor on $\mathcal{C}(\Omega)$. We define an automorphism $T_{\Omega}$
  of $\mathcal{F}$ as $T_{\Omega}=T_{k_n}\cdots T_{k_2}T_{k_1}.$ By
   Theorem 4.7 (4)., $T_{\Omega}$ induces an automorphism of cluster algebra
  $\mathcal{A}_{\Omega}$.

\medskip

\textbf{Definition 4.10.} $T_{\Omega}$ and its inverse
$T^{-}_{\Omega}$ are called the Coxeter automorphisms of the cluster
algebra  $\mathcal{A}_{\Omega}$. $T_{\Omega}$ is simply denoted by
$T$.
\medskip

We recall that $\Omega_0$ and $\Omega _0'$ denote the two
orientations of $\G$ such that in any of these two orientations, any
vertex is a sink or source. For such an orientation $\Omega$, we
denote
 by $\G_{+}$ the set of sinks in $\Omega$, by $\G_{-}$ the set of
sources in $\Omega$. Then $\G=\G_+\cup\G_-.$ Dually we have
$\G=\G_+'\cup\G_-'.$ Now we set $T_{\varepsilon}=\prod _{i\in
\G_{\varepsilon}}T_i$ for $\varepsilon\in \{ +\ -\}.$ Note that
$T=T_+T_-$.
\medskip

\textbf{Corollary 4.11.} $\chi _{\Omega_0}=\chi _{\Omega_0}'; \
T_{\pm}$ are automorphisms of $\mathcal{A}_{\Omega_o}$ and induce
 a bijection from $\chi _{\Omega_0}$ to itself, which sends clusters
to clusters.
\medskip

\textbf{Proof.} We prove firstly that $\chi _{\Omega_0}=\chi
_{\Omega_0'}.$ Let $B_{\Omega_0}$ and $B_{\Omega _0'}$ be the
skew-symmetrizable integer matrices corresponding to the quivers
$(\G, \Omega_0)$ and
 $(\G, \Omega_0')$ respectively. Then $B_{\Omega _0'}=-B_{\Omega
 _0}.$ From the initial seed $(\underline{u}, B_{\Omega_0})$ and
 the initial seed $(\underline{u}, B_{\Omega_0'})$ respectively,
  the new seeds $(\underline{u_k'}, B_{\Omega_0}')$ and
  $(\underline{u_k'}, B_{\Omega_0'}')$ obtained by one step mutation
 in any direction $k$ contain the same
 cluster variables $u_1, \cdots, u_{k-1}, u_k', \cdots, u_n$, and their matrices
 $B_{\Omega_0}'$ and $B_{\Omega_0'}'$ also satisfy the relation
 $B_{\Omega_0'}'=-B_{\Omega_0}'.$
 By induction, we have that $\chi _{\Omega_0}=\chi _{\Omega_0'}.$
  The second and the third statements follow easily from the first one and statements (2) and (4) in
  Theorem 4.7. The proof is finished.

 \medskip

\textbf{Remark. 4.12.} If the number of vertices of $\G$ is $2$, the
automorphisms $T_{\pm}$ of $\mathcal{A}_{\Omega_o}$ are defined in
[SZ], these automorphisms are used to study the positivity and
canonical bases in rank $2$ cluster algebras there.

\medskip

In the following, we will generalize Theorem 4.7. to arbitrary
valued graph $\G$. Let $\Omega$ be an orientation of $\G$,
$C=R^+_{k_n}\cdots R^+_{k_2}R^+_{k_1}$ the corresponding Coxeter
 functor on $\mathcal{C}(\Omega)$ and  $T=T_{k_n}\cdots T_{k_2}T_{k_1}$
 the corresponding Coxeter automorphism of $\mathcal{A}_{\Omega}.$

  We define:
  $$\chi'_{\Omega}=\{ T^m(u_k)\  |  \  m\in \mathbf{Z}, 1\le k \le
  n \  \}.$$

  \medskip

  When $\G$ is of finite type, then $\chi'_{\Omega}=\chi_{\Omega}$
  is the set of cluster variables of $(\G, \Omega)$ by Corollary 3.4.
  and Theorem 4.7.  For infinite type, we prove that elements
  in $\chi'_{\Omega}$ are
  cluster variables of the initial seed $(\underline{u},B)$ and prove the map $\phi_{\Omega}: C^k(P_i)\mapsto T^k(u_i) ,$ for any $k\in\mathbf{Z}, \ i=1,\cdots, n, $  sends tilting objects in $\mathcal{PI}(\Omega)$ to clusters in $\chi'_{\Omega}$ (compare Remark 4.8.).

  \medskip

\textbf{Theorem 4.13.} Let $\G$ be any valued graph, $\Omega$ an
orientation of $\G$. Then  any element in $\chi'_{\Omega}$ is a
cluster variable of the initial seed $(\underline{u},B)$, where $B$ is the
skew-symmetrizable matrix corresponding to $(\G,\Omega)$. Furthermore, the assignment $\phi_{\Omega}:
C^m(P_k[1])\mapsto T^m(u_k), \forall m\in \mathbf{Z}, \  k\in \G $ is a surjection from $\mathcal{PI}(\Omega)$ to
$\chi'_{\Omega}$ such that under this correspondence, $P_k[1]$
corresponds to $u_k$.
\medskip

\textbf{Proof.} Firstly, we show the map $\phi_{\Omega}:
\mathcal{PI}(\Omega)\rightarrow \chi'_{\Omega}$ is well-defined.

 Since any indecomposable object in $\mathcal{C}(\Omega)$ is of the form $P_k[1]$ or $X$ for
 some indecomposable $H-$module $X$ [BMRRT], it follows from [G] [BrB] that $C^m(X)\cong\tau^mX$ for any $X$ in
 ind$\mathcal{C}(\Omega).$

If $\G$ is of finite type,
 then $\chi'_{\Omega}=\chi_{\Omega},$ \
$\mathcal{PI}(\Omega)=\mbox{ind}\mathcal{C}(\Omega)$, and the map
$\phi_{\Omega}$ is well defined, bijective and sends titling objects to clusters by Theorem 4.7.

If $\G$ is of infinite type, then for any pair of indecomposable
objects $X, Y$ in $\mathcal{PI}(\Omega)$, $C^s X\cong C^t Y$ if and
only if $\tau ^sX\cong \tau^tY$ if and only if $s=t$ and $X\cong Y.$
Therefore the map $\phi_{\Omega}$ is well-defined and is surjective,
and under this map, $P_k[1]$
 is sent to $u_k$ for all $k$.

  In the following, we suppose that $\G$
is of infinite type. We will prove that $(T^m(u_1), \cdots,
T^m(u_n))$ is a cluster. For $m=0$, we know $(u_1, \cdots, u_n)$ is
the initial cluster, this corresponds to the slice $(P_1[1],\cdots,
P_n[1])$ in the AR-quiver of the cluster category
$\mathcal{C}(\Omega)$. Since $C^m$ is a triangle auto-equivalence of
the cluster category $\mathcal{C}(\Omega),$ it sends cluster-tilting
set $(P_1[1],\cdots, P_n[1])$ to cluster-tilting set
$(C^mP_1[1],\cdots, C^mP_n[1]).$ By the proof of Theorem 4.7.(2) and
Remark 4.9., we have that $T^m$ sends the initial cluster
$(\underline{u}, B_{\Omega})$ to the cluster $(T^m(\underline{u}),
B_{\Omega}),$ where $T^m(\underline{u})=(T^m(u_1), \cdots,
T^m(u_n))$ is the image of
 $(C^mP_1[1],\cdots, C^mP_n[1])$ under $\phi _{\Omega}.$
 This proves that all elements of the form
$T^m(u_i)$ are cluster variables.
 The proof is finished.
 \medskip

Note that $\phi_{\Omega}$ should be injective, but up to now, we
could not find a proof for this (compare [Z3]).

\medskip

 We remind that $\Phi'_{\ge-1}$ denotes the subset of $\Phi_{\ge-1}$ consisting
of $-\alpha_i, \   i\in \G $ and  \textbf{dim}$X$ for any $X\in
\mathcal{P} \bigcup \mathcal{I},$ where $\mathcal{P}$
 and $\mathcal{I}$ are those in Theorem 3.3.. By Proposition 2.1., $\gamma _{\Omega}$ gives a bijection
from $\mathcal{PI}(\Omega)$ to $\Phi'_{\ge-1}$. Set
$\sigma=\sigma_{i_n}\cdots \sigma_{i_1}.$ It is a bijection of
$\Phi_{\ge -1}.$ Let $P_{\Omega}$ be a map from $\Phi'_{\ge-1}$ to
$\chi'_{\Omega}$ defined as $\phi_{\Omega}\gamma_{\Omega}^{-1}$.
\medskip

\textbf{Proposition 4.14.} Let $\G$ be any valued graph, $\Omega$ an
orientation of $\G$. Then  $P_{\Omega}: \Phi'_{\ge-1} \rightarrow
\chi'_{\Omega}$ is surjective, and the following diagram commutes:
\[ \begin{CD}
\mathcal{PI}(\Omega) @>C
>> \mathcal{PI}(\Omega)\\
@V\gamma_{\Omega} VV  @VV\gamma_{\Omega}V  \\
\Phi'_{\Omega} @>\sigma >> \Phi'_{\Omega}\\
@VP_{\Omega} VV  @VVP_{\Omega}V  \\
\chi'_{\Omega} @>T>> \chi'_{\Omega}
\end{CD}. \]

\medskip

\textbf{Proof.} Since $\phi_{\Omega}$ is surjective and
$\gamma_{\Omega}$ is bijective,
$P_{\Omega}=\phi_{\Omega}\gamma_{\Omega}^{-1}$ is surjective.  The
upper square in the diagram is commutative due to Proposition 3.6..
We verify the commutativity of the lower square.  we verify that
$T\phi_{\Omega}=\phi_{\Omega}C :$  For any $C^m(P_k[1])\in
\mathcal{PI}(\Omega),\
T\phi_{\Omega}(C^m(P_k[1]))=T(T^m(u_k))=T^{m+1}(u_k)=\phi_{\Omega}(C^{m+1}(P_k[1]))
=\phi_{\Omega}C(C^{m}(P_k[1])).$ Then
$TP_{\Omega}=T\phi_{\Omega}\gamma_{\Omega}^{-1}=
\phi_{\Omega}C\gamma_{\Omega}^{-1}=
\phi_{\Omega}\gamma_{\Omega}^{-1}\sigma= P_{\Omega}\sigma.$  The
proof is finished.

\medskip

Let $<C>$ and $<T>$ be the groups generated by $C$ and $T$
respectively. Let $h$ be the Coxeter number of a finite root system
$\Phi$ [FZ3].
\medskip

\textbf{Theorem 4.15.} Let $\G$ be a valued graph of finite type and
$\Omega$ an orientation of $\Delta$. Then the order of $C$ is equal
to $\frac{h+2}{2}$ if the longest element of Weyl group is $-1$, and
is equal to $h+2$ otherwise.
\medskip

\textbf{Proof.} For any two orientations $\Omega, \ \Omega'$, the corresponding Coxeter functors $C_{\Omega}$ and
$C_{\Omega'}$ are conjugated with each other by Proposition 4.14.
  So we choose
$\Omega$ is the orientation such that any vertex is a sink or
source. It follows from Theorem 2.6 in [FZ3] that the order of
$\sigma $ is equal to $\frac{h+2}{2}$ if the longest element of Weyl
group is $-1$, and is equal to $h+2$ otherwise. Applying Proposition
4.14 to $(\G , \Omega )$, we have the order of $C$ is the same as
that of $\sigma$. The proof is finished.

 \medskip

 \textbf{Corollary 4.16.} Let $\G$ be a valued graph of finite type and
$\Omega$ an orientation of $\G$. Then the order of $T$ is equal to
$\frac{h+2}{2}$ if the longest element of Weyl group is $-1$, and is
equal to $h+2$ otherwise.
\medskip

\textbf{Proof.} If $\G$ is of finite type, the maps $\phi_{\Omega}$
and $P_{\Omega}$ in Proposition 4.14. are bijections. Then the order of
$T$ is the same as $C$. This finishes the proof.

\medskip

\textbf{Corollary 4.17.} Let $\G$ be a valued graph of finite type
and $\Omega$ an orientation of $\G$. Then $\chi _{\Omega}=\{ T^{m}(
u_k)\ |\ k\in \G,\ 0\le m\le\frac{h+2}{2} \}$ if the longest element
of Weyl group is $-1$, and $\chi _{\Omega}=\{ T^{m}( u_k)\ |\ k\in
\G,\ 0\le m\le h+2 \}$  otherwise.
\medskip

\textbf{Proof.} This is a consequence of Theorem 4.13. and the
Corollary 4.16.

\medskip

\begin{center}
\textbf {ACKNOWLEDGMENTS.}\end{center} The author would like to
thank Professor Idun Reiten for her helpful conservation on this
topic. The author would like to thank Henning Krause for his
hospitality when he was visiting Paderborn. The author is grateful
to the
 referees for a number of helpful comments and valuable suggestions.
%\newpage

\begin{center}

\end{center}

\medskip
\end{document}